\documentclass{article}

\def\1ox{{ \Omega^1_{\scriptstyle{X}} }}
\def\2ox{{ \Omega^2_{\scriptstyle{X}} }}

\def\ok1{{ \Omega^1_K }}
\def\ok2{{ \Omega^2_K }}
\def\Om{{ \Omega }}

\def\O{{ {\mathcal O} }}

\def\ra{{ \rightarrow }}
\def\da{{ \downarrow }}

\def\a{{ \alpha }}

\def\hra{{ \hookrightarrow }}
\def\da{{ \downarrow }}

\def\F{{ {\bf F} }}

\def\D{{ \Delta }}

\def\8{{ {\infty } }}

\def\^{{ ^{\wedge} }}

\def\K{{ \bar{K} }}

\newtheorem{thm}{Theorem}
\newtheorem{cor}{Corollary}

\usepackage{amsmath}
\usepackage{amsfonts}
\usepackage{amscd}
\usepackage{amssymb}

\def\D{{ \Delta }}

\def\D{{ \Delta }}

\def\K{ \Om^d }
\def\F{ {\bf F} }
\title{A  vanishing theorem for Fano varieties in positive characteristic}
\author{ Minhyong Kim }
\begin{document}
\maketitle
\section{Introduction}
Let $k$ be a perfect field of characteristic $p>0$.
$W=W(k)$ is the ring of Witt vectors of $k$ and $K$ is its
fraction field. 

By  a {\em Fano variety} or an
{\em anti-canonical variety} over $k$, we mean a smooth projective
variety $X$ over $k$ such that the anti-canonical sheaf
$(\K_X)^*$ is ample, where we denote by $d$ the dimension of $X$.

In this note, we wish to prove the following

\begin{thm}
Let $X$ be a Fano variety over $k$ of dimension $d$. Then
$$H^i(X, W\O_X)\otimes K=0$$
for $i>0$.
\end{thm}

By Poincar\'e duality, this implies that
$$H^i(X, W\K_X)\otimes K=0$$
for $i<d$, where
$W\K_X$ is the sheaf of De Rham-Witt (DRW) differential forms
of degree $d$ constructed by Bloch and Illusie. 
This is because $H^i(W\O_X)\otimes K$ is the
part of crystalline cohomology with Frobenius slopes in the interval
$[0,1[$, while $H^i(W, W\K_X)$ is the part with slope $d$
(\cite{Il} cor. II.3.5).
The author gave an erroneous proof of this corollary
(in fact, an integral version) in an earlier preprint.
H. Esnault \cite{Es} subsequently gave a correct proof using the
Bloch-Srinivas decomposition theorem in rational
Chow groups and rigid cohomology.

Here, we will use the additional structure provided
by the De Rham-Witt complex 
(\cite{Bl}, \cite{Il}, \cite{Hyo1},\cite{HK}) in order to prove
the theorem which is a slight strengthening of Esnault's
result. 

In fact, the theorem has the following corollary:

\begin{cor}
Let $X$ be a Fano variety over a finite field
with $q$ elements. Then
the number of rational points on $X$
is congruent to 1 mod $q$.
\end{cor}
The corollary is an easy consequence
of the Lefschetz trace formula for crystalline cohomology
and slope arguments. Esnault's theorem gives that
the number is congruent to 1 mod $p$, if $q=p^n$.
\section{Proof}
Throughout, if we write $H^i$ without further embellishments,
we mean rational crystalline cohomology.

We start with a quick summary of Esnault's
proof (which was an adaptation of Bloch's proof in
characteristic zero \cite{Bl2}): 
Because $X$ is rationally connected and therefore
$CH_0(X)\otimes Q=0$ (\cite{Ko}), one gets from the Bloch-Srinivas
theorem \cite{BlSr}
that the diagonal correspondence $\D \subset X\times X$
is equivalent in $CH_0(X)\otimes Q $ to a sum
$$z\times X + Z$$
where $z$ is the class of a closed point and $Z$ is a cycle supported on
$X\times U$ for $U\subset X$ the complement of
a divisor $D\subset X$.
Therefore, if we apply the diagonal correspondence to
a class of $H^i(X)$, $i>0$, then the only thing that
acts is the $Z$ part. That is, if $i>0$, and $\a \in H^i(X)$ then
$[\D]_*(\a)=[Z]_*(\a)$. We have
$$H^i(X) \simeq H^i_{rig}(X)$$
where $H^i_{rig}$ is Berthelot's rigid cohomology \cite{Be},
and rigid cohomology has nice properties for
open varieties.
So if regard $\a$ as a class in $H^i_{rig}(X)$ and
 pull back to the
subset $U$, then
$\a_U=[Z_{X\times U}]_*(a)=0$ since $Z$ is supported on
$X\times D$. On the other hand, one argues that the map 
$H^i_{rig}(X) \ra H^i_{rig}(U)$
is injective on the Frobenius slope zero part.
This concludes the argument. 

Now we modify this proof. We wish to show
that $H^i_{rig}(X) \ra H^i_{rig}(U)$ is in fact
injective on the part with slope in $[0,1[$.
First, let $f:Y\ra X$ be an alteration \cite{DeJ}
with the property that $E=f^*(D)$ is of normal crossing.
We have a commutative diagram
$$\begin{array}{ccc}
V &\hra & Y\\
\downarrow & & \downarrow\\
U &\hra &X 
\end{array}$$
where $V=Y-E$. Since the pullback from $X$ to $Y$ is
injective, we just need to show that the map
$H^i_{rig}(Y) \ra H^i_{rig}(V)$ is injective on
the part with slope in $[0,1[$. But according to
Shiho's comparison theorem \cite{Sh},
we have a commutative diagram
$$\begin{array}{ccc}
H^i_{rig}(Y)& \ra & H^i_{rig}(V) \\
\da {\scriptstyle \simeq} & & \da {\scriptstyle \simeq} \\
H^i(Y) & \ra & H^i(Y,E)
\end{array}
$$
where the space $H^i(Y,E)$ refers to the rational log crystalline
cohomology of the log scheme $(Y,E)$.
The map
$H^i(Y) \ra  H^i(Y,E)$
is induced by a map
$$W\Om_Y \ra W\Om_Y (\log E)$$
of De Rham-Witt complexes.
This is because we can realize the map at the level of
crystalline complexes \cite{HK} for the two log schemes $Y$ (with trivial
log structure) and $(Y,E)$
and the De Rham-Witt complexes are just given level by
level as the cohomology sheaves of the crystalline complexes.
Meanwhile, the degree zero part is the same and equal to
$W\O_Y$ for both complexes.
So we are done if
we can show that the slope spectral sequence
for $ H^i(Y,E)$ degenerates at $E_1$, as in the case
without log structures, and induces an isomorphism
between $H^i(W\Om^j_Y(\log E))\otimes K$ and the part of
$ H^i(Y,E)$ with slope in $[j,j+1[$.
To see this, one needs only repeat verbatim 
Bloch's argument from \cite{Bl}, III.3.
This is because the log de Rham-Witt complex
$W\Om_Y (\log E)$
is also equipped with operators $V$ and $F$ satisfying
$$FV=VF=p$$
$$pFd=dF, \ \ Vd=pdV$$
which is all that is necessary for Bloch's argument to apply:
In brief, the map given by $p^jF$ on $W\Om^j_Y(\log E)$
is a map of complexes, and induces on log crystalline
cohomology the action of the absolute Frobenius $\phi$.
So on each $E_r$ of the spectral sequence,
the map induced by $p^jF$ on the subquotient
$E^{ji}_r$ of  $H^i(W\Om^j_Y(\log E))\otimes K$
is the map that commutes with the differentials.
So the slope of $\phi$ on $E^{ji}_r$ is $\geq j$.
However, from $FV=p$ and the fact that $V$ acts
topologically nilpotently, we deduce that $\phi|E^{ji}_r$
has slope $<j+1$. Therefore, the difference of slopes
forces all the differentials to be zero from $E_1$ on,
and the Dieudonne-Manin classification of crystals
allows us to split the filtration of the
spectral sequence.

Let us dispense of the easy corollary 1:
As already mentioned,
  $H^i(W\O_X)\otimes K$ is identified
with the part of crystalline cohomology $H^i(X)$
on which the operator
$\phi$ induced by the absolute Frobenius of $X$
has  slope in $[0,1[$. Thus, the vanishing shows that
 all the Frobenius slopes on
$H^i(X)$ are $\geq 1$ for $i>0$. 
Now when
$q=p^n$ and $k$ is the finite field $\F_q$, the Lefschetz trace formula
gives us
$$|X(\F_q)|=\Sigma_i (-1)^i \mbox{Tr} (\phi^n|H^i_{cr}(X)\otimes K)$$
which obviously yields our congruence.

\medskip
{\bf Acknowledgements:} The author was supported in part
by a grant from the NSF. He also benefitted from the 
excellent environment of the MPIM during the
preparation of this note.
As with several previous papers, the author is grateful to
Arthur Ogus for a prompt
response to an email inquiry.

{\footnotesize DEPARTMENT OF MATHEMATICS, UNIVERSITY OF ARIZONA,
TUCSON, AZ 85721, U.S.A. EMAIL: kim@math.arizona.edu}

\end{document}